\documentclass[11pt,reqno]{amsart}
\usepackage[margin=1in]{geometry}
\usepackage{amssymb,amsthm,amsmath,amsfonts,thmtools,bbm}
\usepackage{hyperref}
\usepackage[capitalize]{cleveref}

\newtheorem{theorem}{Theorem}[section]

\newtheorem{lemma}[theorem]{Lemma}
\newtheorem{conjecture}[theorem]{Conjecture}

\theoremstyle{definition}

\newtheorem*{remark}{Remark}

\numberwithin{equation}{section}

\newcommand{\bz}{\mathbf 0}
\newcommand{\CC}{\mathbb C}
\newcommand{\ZZ}{\mathbb Z}
\newcommand{\cV}{\mathcal V}
\newcommand{\1}{\mathbbm 1}
\DeclareMathOperator{\wt}{wt}

\usepackage{xcolor}

\title{Modular $n$-queens in higher dimensions}

\begin{document}

\author{Zachary Halberstam}
\email{zhalberstam@g.harvard.edu}
\address{Department of Economics, Harvard University, Cambridge, MA 02138, USA}
\author{Mitchell Lee}
\author{Carl Schildkraut}
\email{carlsch@stanford.edu}
\address{Department of Mathematics, Stanford University, Stanford, CA 94305, USA}

\thanks{ZH was supported by a Herchel Smith Harvard Scholarship during part of the research. CS is supported by the National Science Foundation Graduate Research Fellowship Program under Grant No.\ DGE-2146755.}

\begin{abstract}
    A chess queen in $d$-dimensions can move any distance in the direction of any edge or diagonal of the $d$-dimensional hypercube.
    (That is, a queen can move from a position $x$ to $x+cv$ for any $c$ and any nonzero $v\in\{-1,0,1\}^d$.)
    We prove that one can place $n^{d-1}$ non-attacking queens on the $d$-dimensional modular board $(\mathbb Z/n\mathbb Z)^d$ if and only if all prime factors of $n$ exceed $2^d$.
    The $d=3$ case of this result is a conjecture of van Rees from 1981.
    Our proof involves both human ideas and substantial AI assistance.
\end{abstract}

\maketitle

\section{Introduction}\label{sec:intro}

The $n$-queens problem is one of the oldest in combinatorics.
The modular version, where attacks may wrap around the edges of the board, dates back at least to 1918, when P\'olya \cite{polya} proved that $n$ mutually non-attacking queens can be placed on a modular $n$-by-$n$ chessboard if and only if $n$ is not a multiple of $2$ or $3$. In this article, we establish an analogous result for the higher-dimensional version of the problem. 

Let $n,d>1$ be integers.
We will consider placing ``queens'' on elements of $(\ZZ/n\ZZ)^d$. Two queens in positions $x,y\in(\ZZ/n\ZZ)^d$ are said to \emph{attack} each other if their difference $x-y$ is a scalar multiple of some nonzero element of $\{-1,0,1\}^d$.
It is clear that one cannot place more than $n^{d-1}$ non-attacking queens in $(\ZZ/n\ZZ)^d$, since no two non-attacking queens may share their first $d-1$ coordinates.
Our main result is to provide necessary and sufficient conditions under which this bound can be attained.

\begin{theorem}\label{thm:main}
    Let $d\geq 2$.
    One can place $n^{d-1}$ non-attacking queens on $(\ZZ/n\ZZ)^d$ if and only if all prime factors of $n$ exceed $2^d$.
\end{theorem}

\subsection{History} In 1980, Klarner proved the sufficiency direction of \cref{thm:main} in the case $d = 3$ \cite{Klarner} with an explicit construction. Klarner's construction naturally generalizes to prove the sufficiency direction of \cref{thm:main} for all $d$ (see, for example, \cite[Corollary~8.2]{HS} and \cite[Theorem~1]{Williams}). The new contribution here is the necessity direction.

The $d=3$ case of \cref{thm:main} was conjectured by van Rees in 1981 \cite[Conjecture-3.2]{vanRees}. Independently, Allison, Yee, and McGaughey made the same conjecture in \cite{aym} (though in the literature the conjecture is usually attributed to van Rees). In their literature review, Bell and Stevens reproduced the conjecture as \cite[Conjecture 24]{bellstevens} and also posed the conjecture for $d\geq 3$ (which \cref{thm:main} verifies) as \cite[Conjecture 25]{bellstevens}.

P\'olya's result \cite{polya} for $d=2$ implies an analogous result for $d>2$: if $n^{d-1}$ non-attacking queens can be placed in $(\ZZ/n\ZZ)^3$, then $2\nmid n$ and $3\nmid n$.
The only other progress in this direction (to our knowledge) is a 2026 result of Williams \cite[Theorem~4]{Williams}, who showed that, for $d=3$, either $5\nmid n$ or $25\mid n$.
Williams' approach is similar to that of P\'olya, in that it involves summing polynomials in $\ZZ[x,y,z]$ over the queen positions $(x,y,z)\in(\ZZ/n\ZZ)^3$.
Our approach is different; we use multiplicative characters of $(\ZZ/n\ZZ)^d$ (that is, Fourier analysis) rather than polynomials. 
This seems to produce stronger results; it is an intriguing question whether the approach with polynomials can be stated in a Fourier-theoretic framework.

We remark additionally that the $n$-queens problem is quite storied and there have been multiple variants of the problem proposed in higher dimensions, some of which we discuss later. 
See the survey of Bell and Stevens \cite{bellstevens} and the master's thesis of Kunt \cite{kunt} for references.

\subsection{Organization}
We give a proof of the sufficiency direction of \cref{thm:main} (which is not new here) in \cref{sec:sufficiency}.
We prove the necessity direction in \cref{sec:necessity}.
In \cref{sec:variants}, we discuss some variants of the $n$-queens problem and some potential extensions of our work.

Our proof of \cref{thm:main} was produced in collaboration with AI.
We outline this process and describe the sources behind the various ideas in the proof in Appendix~\ref{sec:appendix}.

\section{Sufficiency}\label{sec:sufficiency}

In this section, we prove the sufficiency direction of \cref{thm:main}: if all prime factors of $n$ are greater than $2^d$, then $n^{d-1}$ non-attacking queens can be placed on $(\ZZ/n\ZZ)^d$.
This argument is a simple extension of the arguments used in the $d=2$ case and appears, for example, as \cite[Theorem~1]{Williams} and \cite[Corollary~8.2]{HS}.
The argument in the $d=3$ case is due to Klarner \cite{Klarner}.
We present the full argument here for completeness.

\begin{proof}[Proof of \cref{thm:main}, sufficiency direction]
	Let $\xi:=(1,2,4,\ldots,2^{d-1})$.
	Place a queen in position $x \in (\ZZ/n\ZZ)^d$ if and only if $\langle\xi,x\rangle=0$. It is clear that $n^{d - 1}$ queens are placed in this way.
	Assume for the sake of contradiction that two queens at positions $x$ and $y$ are attacking. Then, there exists some nonzero $j\in\ZZ/n\ZZ$ and $v\in\{-1,0,1\}^d\setminus\{\bz\}$ for which $x-y=jv$. 
	We conclude that
	\[0=\langle\xi,x\rangle-\langle\xi,y\rangle=j\langle\xi,v\rangle=j\sum_{k=0}^{d-1}2^kv_k.\]
	Let $r=\sum_{k=0}^{d-1}2^kv_k$. Write $v=u-w$ for distinct $u,w\in\{0,1\}^d$. 
	By the uniqueness of binary representation, we have $\sum_{k=0}^{d-1}2^ku_k\neq\sum_{k=0}^{d-1}2^kw_k$. We conclude that, as an integer, $r$ is nonzero.
	Moreover, $\lvert r\vert\leq 2^d-1$.
	Therefore, by our assumption on $n$, we have $\gcd(r,n)=1$. The property $jr=0$ for $j\in\ZZ/n\ZZ$ nonzero provides a contradiction.
\end{proof}

\begin{remark}
    The configuration of $n^{d-1}$ non-attacking queens presented in this proof is an \emph{affine} configuration; that is, the cells at which we place queens form a hyperplane $\{x:\langle\xi,x\rangle=0\}$ in $(\ZZ/n\ZZ)^d$.
    This proof can be easily extended to show that, in fact, such an affine configuration of $n^{d-1}$ non-attacking queens only exists if $n$ has no prime factors less than $2^d$.

    Indeed, let $\xi\in(\ZZ/n\ZZ)^d$, and suppose that $p<2^d$ is a prime factor of $n$.
    Since $p<2^d$, there exist distinct $(u_1,\ldots,u_d),(w_1,\ldots,w_d)\in\{0,1\}^d$ such that
    \[\sum_{k=1}^du_k\xi_k\equiv \sum_{k=1}^dw_k\xi_k\pmod p.\]
    This implies that
    \[\left\langle\xi,\left(\frac np(u_1-w_1),\ldots,\frac np(u_d-w_d)\right)\right\rangle=0.\]
    So, the set $\{x:\langle\xi,x\rangle=0\}$ contains the attacking positions $(0,\ldots,0)$ and $(n(u_1-w_1)/p,\ldots,n(u_d-w_d)/p)$, and is thus not a set of non-attacking queens.

    Given this, \cref{thm:main} can be rephrased in the following way: given $n$ and $d$, there exists a configuration of $n^{d-1}$ non-attacking queens on $(\ZZ/n\ZZ)^d$ if and only if there exists such an affine configuration.
\end{remark}

\section{Necessity}\label{sec:necessity}

We now verify the other direction of \cref{thm:main}.
We first set some notation which we will carry through the proof.
Let $\Gamma:=(\ZZ/n\ZZ)^d$ and let $\cV_d:=\{-1,0,1\}^d\setminus\{\bz\}$. Write $\omega:=e^{2\pi i/n}$.
For $f\colon\Gamma\to\CC$, define the Fourier transform $\widehat f \colon \Gamma \to \CC$ by
\[\widehat f(\xi):=\sum_{x\in\Gamma}f(x)\omega^{\langle\xi,x\rangle}.\]
Define the convolution of two functions $f, g \colon \Gamma \to \CC$ by
\[(f*g)(x)=\sum_{y\in\Gamma}f(y)g(x-y).\]
We have $\widehat{(f*g)}(\xi)=\widehat f(\xi)\widehat g(\xi)$.
For $y\in\Gamma$, write $\delta_y$ for the function which evaluates to $1$ on $y$ and $0$ everywhere else. Note that $\widehat{\delta_y}(\xi)=\omega^{\langle\xi,y\rangle}$.

If $S \subseteq \Gamma$ is a set of positions of $n^{d-1}$ non-attacking queens, then each line $\{x+jv:j\in\ZZ/n\ZZ\}$ for $x\in\Gamma$ and $v\in\cV_d$ contains exactly one element of $S$.
We will use the following lemma, which is a special case of \cite[Lemma~7.4]{HS}, to conclude some important properties of the Fourier transform of the indicator function of $S$.

\begin{lemma}\label{lem:fhat}
	Let $f\colon\Gamma\to\ZZ$.
    Let $C\in\ZZ$.
	Suppose that, for every $x\in\Gamma$ and every $v\in\cV_d$, we have
	\[\sum_{j\in\ZZ/n\ZZ}f(x+jv)=C.\]
    Let $\xi\in\widehat\Gamma\setminus\{0\}$. 
    If $\langle\xi,v\rangle=0$ for some $v\in\cV_d$, then $\widehat f(\xi)=0$.
\end{lemma}
\begin{proof}
    For each $x\in\Gamma$, we have
    \[\sum_{j\in\ZZ/n\ZZ}f(x+jv)\omega^{\langle\xi,x+jv\rangle}=\sum_{j\in\ZZ/n\ZZ}f(x+jv)\omega^{\langle\xi,x\rangle}=C\omega^{\langle\xi,x\rangle}.\]
    Therefore
    \[\widehat f(\xi)=\sum_{x\in\Gamma}f(x)\omega^{\langle\xi,x\rangle}=\frac1n\sum_{x\in\Gamma}\sum_{j\in\ZZ/n\ZZ}f(x+jv)\omega^{\langle\xi,x+jv\rangle}=\frac Cn\sum_{x\in\Gamma}\omega^{\langle\xi,x\rangle}=0.\qedhere\]
\end{proof}

We now explain the strategy behind the proof of \cref{thm:main}.
Let $S \subseteq \Gamma$ be a set of $n^{d-1}$ non-attacking positions.
\cref{lem:fhat} implies that, if $\xi\in\widehat\Gamma\setminus\{0\}$ is orthogonal to some element of $\cV_d$, then $\widehat{\1_S}(\xi)=0$.
If $n<2^d$, then this condition implies that $\widehat{\1_S}(\xi)=0$ for every nonzero $\xi\in\Gamma$. Thus $\1_S$ is constant, which cannot be.

We will generalize this argument by constructing a function $g\colon\Gamma\to\ZZ$ with the property that $\widehat g(\xi)=0$ whenever $\xi\in\widehat\Gamma\setminus\{0\}$ is not orthogonal to any element of $\cV_d$. 
For any such $g$, by \cref{lem:fhat}, the function $\1_S*g$ will be constant and integer-valued, and so we will have
\[n^d\mid\sum_{x\in\Gamma}(\1_S*g)(x)=\widehat{(\1_S*g)}(0)=\widehat{\1_S}(0)\widehat g(0)=n^{d-1}\widehat g(0).\]
If $n\nmid\widehat g(0)$, this provides the desired contradiction.

To construct such a $g$, we first prove an auxiliary lemma. 
Given a vector $\alpha\in\ZZ_{\geq 0}^r$ for some $r$, we write $\wt\alpha$ for the sum of coordinates of $\alpha$.

\begin{lemma}\label{lem:poly}
	For a positive integer $r\geq2$, write
	\[h_{n}(z_1,\ldots,z_r):=\sum_{\substack{\alpha\in\ZZ_{\geq 0}^r\\\wt\alpha=n-r+1}}z_1^{\alpha_1}\cdots z_r^{\alpha_r}\] for the $n$th complete homogeneous symmetric polynomial in the variables $z_1, \ldots, z_r$.
	If $z_1,\ldots,z_r$ are distinct $n$th roots of unity, then $h_{n}(z_1,\ldots,z_r)=0$.
\end{lemma}
\begin{proof}
	Let $t_1,\ldots,t_r$ be indeterminates. We have
	\[\sum_{n=r-1}^\infty h_{n}(t_1,\ldots,t_r)=\sum_{\alpha\in\ZZ_{\geq 0}^r}t_1^{\alpha_1}\cdots t_r^{\alpha_r}=\prod_{i=1}^r\frac1{1-t_i}.\]
	Now, define the determinants
	\[
	D_m:=
	\begin{vmatrix}
		1&t_1&t_1^2&\cdots&t_1^{r-2}&t_1^m\\
		1&t_2&t_2^2&\cdots&t_2^{r-2}&t_2^m\\
		\vdots&\vdots&\vdots&\ddots&\vdots&\vdots\\
		1&t_r&t_r^2&\cdots&t_r^{r-2}&t_r^m
	\end{vmatrix}.
	\]
	We compute
	\[\sum_{m=0}^\infty D_m=\begin{vmatrix}
		1&t_1&t_1^2&\cdots&t_1^{r-2}&\frac1{1-t_1}\\
		1&t_2&t_2^2&\cdots&t_2^{r-2}&\frac1{1-t_2}\\
		\vdots&\vdots&\vdots&\ddots&\vdots&\vdots\\
		1&t_r&t_r^2&\cdots&t_r^{r-2}&\frac1{1-t_r}
	\end{vmatrix}=\begin{vmatrix}
	1-t_1&t_1-t_1^2&\cdots&t_1^{r-2}-t_1^{r-1}&1\\
	1-t_2&t_2-t_2^2&\cdots&t_2^{r-2}-t_1^{r-1}&1\\
	\vdots&\vdots&\ddots&\vdots&\vdots\\
	1-t_r&t_r-t_r^2&\cdots&t_r^{r-2}-t_r^{r-1}&1
	\end{vmatrix}\prod_{i=1}^r\frac1{1-t_i}.\]
	Performing column operations and permuting the columns,
	\[\sum_{m=0}^\infty D_m=\begin{vmatrix}1&t_1&\cdots&t_1^{r-1}\\1&t_2&\cdots&t_2^{r-1}\\\vdots&\vdots&\ddots&\vdots\\1&t_r&\cdots&t_r^{r-1}\end{vmatrix}\prod_{i=1}^r\frac1{1-t_i}=D_{r-1}\prod_{i=1}^r\frac1{1-t_i}.\]
	We conclude that
	\[\sum_{n=r-1}^\infty h_{n}(t_1,\ldots,t_r)D_{r-1}=\sum_{m=0}^\infty D_m.\]
	The $n$ term on the left is homogeneous of degree $n-r+1+\binom r2$, while the $m$ term on the right is homogeneous of degree $m+\binom{r-1}2$. 
	We conclude that $h_{n}(t_1,\ldots,t_r)=D_n/D_{r-1}$. (This is a special case of Jacobi's bialternant formula for Schur polynomials.)
	
	Now, evaluate this identity at distinct $n$th roots of unity $z_1,\ldots,z_r$. Since the first and last columns of the matrix defining $D_n$ are equal, we have $D_n=0$. On the other hand, $D_{r-1}$ is the determinant of a Vandermonde matrix, so the distinctness of $z_1,\ldots,z_r$ implies that $D_{r-1}\neq 0$.
	The result follows.	
\end{proof}

We now construct our auxiliary function $g$.

\begin{lemma}\label{lem:ghat}
	Let $r\leq 2^d$, and let $y_1,\ldots,y_r$ be distinct elements of $\{0,1\}^d$. Define a function $g\colon\Gamma\to\ZZ$ by
	\[g=\sum_{\substack{\alpha\in\ZZ_{\geq 0}^r\\\wt\alpha=n-r+1}}\delta_{\alpha_1y_1+\cdots+\alpha_r y_r}.\]
    Let $\xi\in\widehat\Gamma\setminus\{0\}$.
    If $\langle\xi,v\rangle\neq 0$ for every $v\in\cV_d$, then $\widehat g(\xi)=0$.
\end{lemma}
\begin{proof}
    For $1\leq j\leq r$, let $\lambda_j=\omega^{\langle\xi,y_j\rangle}$. 
    If $\lambda_j=\lambda_k$ for some $j\neq k$, then $\langle\xi,y_j-y_k\rangle=0$, which would contradict our assumption that $\xi$ is not orthogonal to any element of $\cV_d$.
 	Therefore, $\lambda_1,\ldots,\lambda_r$ are distinct $n$th roots of unity.
 	We now compute
 	\[
 		\widehat g(\xi)=
 		\sum_{
 			\substack{\alpha\in\ZZ_{\geq 0}^r\\\wt\alpha=n-r+1}
 		} \widehat{\delta}_{\alpha_1y_1+\cdots+\alpha_r y_r}(\xi)=
 		\sum_{
 			\substack{\alpha\in\ZZ_{\geq 0}^r\\\wt\alpha=n-r+1}
 		} \omega^{\langle\xi,\alpha_1y_1+\cdots+\alpha_r y_r\rangle}
 		=h_n(\lambda_1,\ldots,\lambda_r).\]
 	The result thus follows from \cref{lem:poly}.
\end{proof}

\begin{proof}[Proof of \cref{thm:main}, necessity direction]
	Let $S\subseteq\Gamma$ be a set of $n^{d-1}$ non-attacking queens and let $p < 2^d$ be a prime. We will prove that $p \nmid n$.
    
    Since $S$ has size $n^{d-1}$, every line $\{x+jv:j\in\ZZ/n\ZZ\}$ for $v\in\cV_d$ has exactly one element of $S$.
    In particular, the indicator function $\1_S$ of $S$ satisfies the precondition of \cref{lem:fhat}.
	
	Now, let $y_1,\ldots,y_{p+1}$ be distinct elements of $\{0,1\}^d$ (we can choose such elements since $p<2^d$). 
	Let $g$ be as in \cref{lem:ghat}.
	\cref{lem:fhat} implies that $\widehat{\1_S}(\xi)=0$ whenever $\xi\neq 0$ is orthogonal to some element of $\cV_d$, while \cref{lem:ghat} implies that $\widehat g(\xi)=0$ for any $\xi\neq 0$ not orthogonal to any element of $\cV_d$.
	We conclude that $\widehat{\1_S*g}(\xi)= \widehat{\1_S}(\xi) \widehat{g}(\xi) = 0$ for every $\xi\in\widehat\Gamma\setminus\{0\}$.
	So, $\1_S*g$ is constant.
	Since $\1_S$ and $g$ take only integer values, $\1_S*g$ takes only integer values.
	We conclude that $n^d\mid \widehat{\1_S*g}(0)$.
	
	On the other hand,
	\[\widehat{\1_S*g}(0)=\widehat{\1_S}(0)\widehat g(0)=n^{d-1}\cdot\left(\sum_{\substack{\alpha\in\ZZ^{r+1}_{\geq 0}\\\wt\alpha=n-p}}1\right)=\binom npn^{d-1}=\frac1pn^d\binom{n-1}{p-1}.\]
	Hence, $\frac{1}{p} \binom{n-1}{p-1}$ is an integer, which implies $p \nmid n$, as desired.	
\end{proof}

\section{Extensions and variants}\label{sec:variants}

There remain many open problems about configurations of non-attacking queens in $(\ZZ/n\ZZ)^d$. Moreover, there are many other variants of the $n$-queens problem in higher dimensions. We discuss some of those extensions and variants here.

If $n$ has a prime factor less than $2^d$, then \cref{thm:main} implies that at most $n^{d - 1}-1$ non-attacking queens can be placed in $(\ZZ/n\ZZ)^d$. It is natural to ask precisely how many can be placed.
In two dimensions, this question is answered by a classical result of Monsky \cite{monsky}.
Williams has made some partial progress towards this question for $d>2$; see \cite[Theorem~3]{Williams} and also \cite[Conjecture~2]{Williams}.
On the other hand, if $n$ has no prime factor less than $2^d$, then there is another natural question: in how many ways can $n^{d-1}$ non-attacking queens can be placed on $(\ZZ/n\ZZ)^d$? The $2$-dimensional version of case has recently been resolved up to a $(1+o(1))^n$ factor by Bowtell and Keevash \cite{bowtellkeevash}.
Neither of these questions seem particularly amenable to the methods we have introduced here.

Probably the most standard variant of the $n$-queens problem (and certainly the most standard in $d=2$) is the ``non-toroidal'' version of the problem we have considered.
Let $n$ and $d$ be positive integers.
Two queens in positions $x,y\in\{1,2,\ldots,n\}^d$ are said to attack each other if their difference is a scalar multiple of some $\{-1,0,1\}^d\setminus\{\bz\}$. 
(For example, the positions $(n-1,0,0)$ and $(0,1,1)$ are non-attacking for $n>2$, while they are attacking in the toroidal problem on $(\ZZ/n\ZZ)^3$.)
Here it remains clear that no more than $n^{d-1}$ non-attacking queens can be placed.
For $d=3$, there is a conjecture of van Rees \cite[Conjecture 3.3]{vanRees} that the existence problem is equivalent to the corresponding existence problem for a toroidal board.

\begin{conjecture}\label{conj:non-toroidal}
    If $\gcd(n,210)>1$, one cannot place $n^2$ non-attacking queens in $\{1,2,\ldots,n\}^3$.
\end{conjecture}

\noindent We believe that \cref{conj:non-toroidal} is false, and that in fact, for fixed $d\geq 2$, one can place $n^{d-1}$ non-attacking queens in $\{1,2,\ldots,n\}^d$ for all sufficiently large $n$. 
Like the two-dimensional $n$-queens problem, higher-dimensional analogues can be phrased as an $r$-partite $r$-uniform hypergraph matching problem. Set $r=(3^d-1)/2$. Construct a hypergraph $\mathcal H$ as follows: create one vertex of $\mathcal H$ for each line in one of the directions $v\in\{-1,0,1\}^d\setminus\{\bz\}$, and create, for each element $x\in\{1,2,\ldots,n\}^d$, an edge of $\mathcal H$ containing the $r$ lines passing through $x$.
A set of non-attacking queen placements corresponds to a matching in $\mathcal H$, and a set of $n^{d-1}$ queens is a maximum-size matching, which must necessarily cover all vertices corresponding to lines in the $d$ coordinate directions.
Following the significant applications of the technique of absorption to the standard and toroidal $n$-queens problem in two dimensions (see \cite{luriasimkin,bowtellkeevash,simkin}), it seems plausible that, at least for $d=3$, where there are only three classes of vertices (corresponding to the three coordinate directions) which require such a matching to be perfect, that such a statement could be proven using absorption.
(Alternatively, such a statement could be proven for some particular $n$ not relatively prime to $210$ by an explicit example.
We have attempted, without any luck, to find such an example using computational tools.)

Other attacking rules also exist for toroidal $n$-queens problems. 
Nudelman \cite{nudelman} introduced a different variant, wherein queens in cells $x$ and $y$ attack each other if $x$ and $y$ lie in the same hyperplane orthogonal to some $v\in\{-1,0,1\}^d\setminus\{\bz\}$, rather than the same line parallel to such a $v$. (In $d=2$, these both collapse to the standard definition.) 
Nudelman proved that, if $n$ is divisible by no prime less than $2^d$, then $n$ ``hyperplane queens'' can be placed in $(\ZZ/n\ZZ)^d$ \cite[Theorem~4]{nudelman}, but that no more than $n-1$ can be placed if $n$ has a prime factor less than $2d$ \cite[Theorem~11]{nudelman}.
It remains an intriguing open question (and one which may be amenable to our techniques) to determine for which $n$ one can place $n$ non-attacking hyperplane queens in $(\ZZ/n\ZZ)^d$.

\section*{Acknowledgements}

We thank Brett Stevens for answering questions about the literature surrounding $n$-queens problems. We thank Jacob Fox and Anubhab Ghosal for helpful advice.

\bibliographystyle{alpha}
\bibliography{bib.bib}

\appendix
\section{AI Disclosure, project history, and methodology}\label{sec:appendix}

As mentioned in \cref{sec:intro}, the main result of this paper was produced in collaboration with AI tools. 
Here we give a history of the project and describe our AI use as well as our process.
We include this information in accordance with a belief that, as AI becomes a larger part of doing mathematics and as process and result are becoming disentangled, it may be useful to report more about the process that led to a result.

ZH and CS first encountered \cite{bellstevens} and the conjecture \cite[Conjecture-3.2]{vanRees} of van Rees (which the $d=3$ case of \cref{thm:main} verifies) in the process of working on \cite{HS} in the summer of 2023, which led to engagement with the modular $n$-queens problem in $d$ dimensions. Initially we believed this conjecture to be false: the sort of argument used by P\'olya \cite{polya} to resolve the $d=2$ case seemed not to extend to $d>3$, and we saw no reason that the classical sufficiency condition (arising from the construction we presented in \cref{sec:sufficiency}) should also turn out to be necessary.

All three authors worked intermittently on the problem in the following years. In particular, we (without the use of AI) pushed the polynomial-type methods used in \cite{polya} and \cite{Williams}, and were unable to obtain any result stronger than \cite[Theorem 4]{Williams}. Also without AI, and building on our earlier work in \cite{HS}, we obtained a partial result towards \cref{thm:main} via Fourier analysis.
Namely, we proved the following result.

\newtheorem{innercustomthm}{Theorem}
\newenvironment{customthm}[1]
  {\renewcommand\theinnercustomthm{#1}\innercustomthm}
  {\endinnercustomthm}
\begin{customthm}{\ref*{thm:main}'}\label{thm:prime}
    Suppose that the prime factorization of $n$ contains a prime power $q$ which is less than $2^d$.
    Then one cannot place $n^{d-1}$ non-attacking queens on $(\ZZ/n\ZZ)^d$.
\end{customthm}
\noindent Theorem~\ref{thm:prime} follows from applying \cref{lem:fhat} to the function $\tilde f\colon(\ZZ/q\ZZ)^d\to\ZZ$ whose values $\tilde f(x)$ count the number of queens in $(\ZZ/n\ZZ)^d$ at positions equivalent to $x$ modulo $q$, and noting that $\tilde f$ cannot be constant if $q^d\nmid n^{d-1}$.

At this point, we continued to believe a counterexample to the van Rees conjecture was possible and considered the computational search problem of finding a counterexample. 
We believed if there were to be a counterexample, there would likely be one at $n=25$, which is the smallest modulus not covered by Theorem~\ref{thm:prime}. We ran these searches starting in 2025, but increased our efforts in May 2026. 
CS ran searches using AlphaEvolve \cite{alphaevolve}, which produced configurations of $25^2$ queens in $(\ZZ/25\ZZ)^3$ with few, but not zero, attacking pairs. 
We tried a mix of strategies, over roughly 40 runs of 12 hours or more (on individual laptops), but none of our configurations came particularly close to working.

This was powerful evidence for the conjecture, so we tried to prove it. 
First, we attempted to extend our Fourier-analytic proof of Theorem~\ref{thm:prime} to the full \cref{thm:main}, but we were unsuccessful.
We then began to use AI tools for this purpose.
First, we queried Codex, which quickly rederived Theorem~\ref{thm:prime}. Then, Codex derived impossibility proofs for the specific cases of $(n,d)=(25,3)$ and $(n,d)=(49,3)$.
The strategy presented by Codex was as follows.
A configuration of $n^2$ queens in $(\ZZ/n\ZZ)^3$ gives a solution to $Bx = \mathbf 1 \pmod n$ for a ``point-line incidence'' matrix $B\in\ZZ^{13n^2\times n^3}$.
Codex then found a vector $c\in\ZZ^{13n^2}$ such that $c^\intercal B=0\pmod m$ but $c^\intercal\mathbf 1\not\equiv 0\pmod m$ for some power $m$ of $n$. (Call such a $c$ a \emph{certificate}.)
We observed that this computation for $(n,d)=(25,3)$ (resp.\ $(n,d)=(49,3)$) is enough to show \cref{thm:main} in the case of $d=3$ and $\nu_5(n)=2$ (resp.\ $\nu_7(n)=2$).
Next, Codex obtained a certificate for $(n,d)=(125,3)$, which bore some similarity to the certificate for $(25,3)$. Codex pointed out a connection of the certificates to Bockstein homomorphisms in homological algebra, which gave us some confidence that this strategy could be extended to all powers of $5$. 
We asked Codex to probe the relationships between the certificates for powers of $5$; eventually, it obtained two distinct proofs of \cref{thm:main} in the case of $(n,d)=(5^k,3)$. 
Finally, Codex produced a proof of the main result. 

We reinterpreted this proof, explicated which ingredients had appeared in our non-AI-assisted work and which ingredients were new, and wrote up the exposition which appears in this article. 
The exposition is entirely ours and we take full responsibility for all the content.
\end{document}